\documentclass[12pt,a4paper,onecolumn]{article}
\usepackage{textcomp}
\usepackage{tikz}
\usepackage{amssymb}
\usepackage{caption}
\usepackage{color}
\usepackage{verbatim}
\usepackage{fancybox}
\usepackage{titlesec} 
\renewcommand\thesection{\arabic{section}} 
\usepackage{setspace}
\usepackage{mathtools}
\usepackage{amsmath}
\usepackage{mathrsfs}
\usetikzlibrary{positioning,fit,calc}
\usepackage[a4paper,left=1in,top=1in,right=1in,bottom=1in,nohead]{geometry} 
\usepackage[hidelinks]{hyperref}
\usepackage{array}
\usepackage{bm}
\usepackage[symbol]{footmisc} 
\def\correspondingauthor{\footnote{Corresponding author. Email: williewong088@gmail.com.}}

\tikzset{block/.style={draw,thick,text width=2cm,minimum height=1cm,align=center},
         line/.style={-latex}}
\newcolumntype{P}[1]{>{\centering\arraybackslash}m{#1}} 

\titleformat{\section}[block]{\large\scshape\bfseries}{\thesection.}{1em}{} 
\titleformat{\subsection}[block]{\bfseries}{\thesubsection.}{1em}{} 

\newtheorem{defn}{Definition}[section]
\newtheorem{thm}[defn]{Theorem}
\newtheorem{eg}[defn]{Example}
\newtheorem{ppn}[defn]{Proposition}

\newtheorem{lem}[defn]{Lemma}

\newtheorem{prob}[defn]{Problem}

\begin{document}
\pagenumbering{arabic}
\begin{center}
    \textbf{\Large A complete characterisation of \\vertex-multiplications of trees with diameter 5}
\vspace{0.1 in} 
    \\{\large W.H.W. Wong\correspondingauthor{}, E.G. Tay}
\vspace{0.1 in} 
\\National Institute of Education\\Nanyang Technological University, Singapore
\end{center}

\begin{abstract}
\indent\par Koh and Tay \cite{KKM TEG 8} introduced a new family of graphs, $G$ vertex-multiplications, as an extension of complete $n$-partite graphs. They proved a fundamental classification of $G$ vertex-multiplications into three classes $\mathscr{C}_0, \mathscr{C}_1$ and $\mathscr{C}_2$. It was shown in \cite{KKM TEG 11} that any vertex-multiplication of a tree with diameter at least 3 does not belong to the class $\mathscr{C}_2$. Furthermore, for vertex-multiplications of trees with diameter $5$, some necessary and sufficient conditions for $\mathscr{C}_0$ were established. In this paper, we give a complete characterisation of vertex-multiplications of trees with diameter $5$ in $\mathscr{C}_0$ and $\mathscr{C}_1$.
\end{abstract}
\section{Introduction}
\indent\par Let $G$ be a graph with vertex set $V(G)$ and edge set $E(G)$. In this paper, we consider graphs $G$ with no loops nor parallel edges, unless otherwise stated.  For any vertices $v,x\in V(G)$, the $\textit{distance}$ from $v$ to $x$, $d_G(v,x)$, is defined as the length of a shortest path from $v$ to $x$. For $v\in V(G)$, its $\textit{eccentricity}$ $e_G(v)$ is defined as $e_G(v):=\max\{d_G(v,x)|\ x\in V(G)\}$. The $\textit{diameter}$ of $G$, denoted by $d(G)$, is defined as $d(G):=\max\{e_G(v)|\ v\in V(G)\}$ while the \textit{radius} of $G$, denoted by $r(G)$, is defined as $r(G):=\min\{e_G(v)|\ v\in V(G)\}$. The above notions are defined similarly for a digraph $D$. A vertex $x$ is said to be \textit{reachable} from another vertex $v$ if $d_D(v,x)<\infty$. For a digraph $D$, the \textit{outset} and \textit{inset} of a vertex $v\in V(D)$ are defined to be $O_D(v):=\{x\in V(D)|\text{ } v\rightarrow x\}$ and $I_D(v):=\{y\in V(D)|\text{ } y\rightarrow v\}$ respectively. If there is no ambiguity, we shall omit the subscript for the above notations.
\indent\par An $\textit{orientation}$ $D$ of a graph $G$ is a digraph obtained from $G$ by assigning a direction to every edge $e\in E(G)$. An orientation $D$ of $G$ is said to be \textit{strong} if every two vertices in $V(D)$ are mutually reachable. An edge $e\in E(G)$ is a \textit{bridge} if $G-e$ is disconnected. Robbins' well-known One-way Street Theorem  \cite{RHE} states the following.
\begin{thm}(Robbins \cite{RHE})
~\\Let $G$ be a connected graph. Then, $G$ has a strong orientation if and only if $G$ is bridgeless.
\end{thm}
Given a connected and bridgeless graph $G$, let $\mathscr{D}(G)$ be the family of strong orientations of $G$. The $\textit{orientation number}$ of $G$ is defined as 
\begin{align*}
\bar{d}(G):=\min\{d(D)|\ D\in \mathscr{D}(G)\}.
\end{align*}
\indent\par The general problem of finding the orientation number of a connected and bridgeless graph is very difficult. Moreover, Chv{\'a}tal and Thomassen \cite{CV TC} proved that it is NP-hard to determine if a graph admits an orientation of diameter 2. Hence, it is natural to focus on special classes of graphs. The orientation number was evaluated for various classes of graphs, such as the complete graphs \cite{BF TR,MSB,PJ 1} and complete bipartite graphs \cite{GG 1,SL}. 
\indent\par In 2000, Koh and Tay \cite{KKM TEG 8} introduced a new family of graphs, $G$ vertex-multiplications, and extended the results on complete $n$-partite graphs. Let $G$ be a given connected graph with vertex set $V(G)=\{v_1,v_2,\ldots, v_n\}$. For any sequence of $n$ positive integers $(s_i)$, a $G$ vertex-multiplication, denoted by $G(s_1, s_2,\ldots, s_n)$, is the graph with vertex set $V^*=\bigcup_{i=1}^n{V_i}$ and edge set $E^*$, where $V_i$'s are pairwise disjoint sets with $|V_i|=s_i$, for $i=1,2,\ldots,n$, and for any $u,v\in V^*$, $uv\in E^*$ if and only if $u\in V_i$ and $v\in V_j$ for some $i,j\in \{1,2,\ldots, n\}$ with $i\neq j$ such that $v_i v_j\in E(G)$. For instance, if $G\cong K_n$, then the graph $G(s_1, s_2,\ldots, s_n)$ is a complete $n$-partite graph with partite sizes $s_1, s_2,\ldots, s_n$. Also, we say $G$ is a parent graph of $G(s_1, s_2,\ldots, s_n)$.
\indent\par For $i=1,2,\ldots, n$, we denote the $x$-th vertex in $V_i$ by $(x,v_i)$, i.e. $V_i=\{(x,v_i)|\ x=1,2,\ldots,s_i\}$. Hence, two vertices $(x,v_i)$ and $(y,v_j)$ in $V^*$ are adjacent in $G(s_1,s_2,\ldots, s_n)$ if and only if $i\neq j$ and $v_i v_j\in E(G)$. For convenience, we write $G^{(s)}$ in place of $G(s,s,\ldots,s)$ for any positive integer $s$, and it is understood that the number of $s$'s is equal to the order of $G$, $n$. Thus, $G^{(1)}$ is simply the graph $G$ itself.

\indent\par The following theorem by Koh and Tay \cite{KKM TEG 8} provides a fundamental classification on $G$ vertex-multiplications.
\begin{thm} (Koh and Tay \cite{KKM TEG 8}) \label{thmA3.1.2}
\\Let $G$ be a connected graph of order $n\ge 3$. If $s_i\ge 2$ for $i=1,2,\ldots, n$, then $d(G)\le \bar{d}(G(s_1,s_2,\ldots,s_n))\le d(G)+2$.
\end{thm}
\indent\par In view of Theorem \ref{thmA3.1.2}, all graphs of the form $G(s_1,s_2,\ldots, s_n)$, with $s_i\ge 2$ for all $i=1,2,\ldots, n$, can be classified into three classes $\mathscr{C}_j$, where 
\begin{align*}
\mathscr{C}_j=\{G(s_1,s_2,\ldots, s_n)|\text{ }\bar{d}(G(s_1,s_2,\ldots, s_n))=d(G)+j\},
\end{align*}
for $j=0,1,2$. Henceforth, we assume $s_i\ge 2$ for $i=1,2,\ldots, n$, unless otherwise stated. The following lemma was found useful in proving Theorem \ref{thmA3.1.2}.

\begin{lem}(Koh and Tay \cite{KKM TEG 8}) \label{thmA3.1.3}
\\Let $\mu_i,\lambda_i$ be integers such that $\mu_i\le\lambda_i$ for $i=1,2,\ldots, n$. If the graph $G(\mu_1,\mu_2,\ldots, \mu_n)$ admits an orientation $F$ in which every vertex $v$ lies on a cycle of length not exceeding $m$, then $\bar{d}(G(\lambda_1, \lambda_2,\ldots, \lambda_n))\le \max\{m, d(F)\}$.
\end{lem}
\indent\par Koh and Tay \cite{KKM TEG 11} further investigated  vertex-multiplications of trees. Since trees with diameter at most 2 are parent graphs of complete bipartite graphs, which are completely solved, they considered trees of diameter at least 3. On the other hand, Ng and Koh \cite{NKL KKM} examined vertex-multiplications of cycles. To continue the discussion, we need some notations.
\indent\par From here onwards, we consider a tree $T$ of diameter $5$ with vertex set $V(T)=\{v_1,v_2,\ldots, v_n\}$. We let $v_1$ and $v_2$ be the two central vertices of $T$, i.e. $e_T(v_k)=rad(T)=3$ for $k=1,2$. For $k=1,2$, denote the neighbours of $v_k$, excluding $v_{3-k}$, by $[i]_k$. i.e. $N_T(v_k)-\{v_{3-k}\}=\{[i]_k|\ i=1,2\ldots, deg_T(v_k)-1 \}$. For each $i =1,2,\ldots, deg_T(v_k)-1$, we denote the neighbours of $[i]_k$, excluding $v_k$ itself, by $[\alpha, i]_k$. i.e. $N_T([i]_k)-\{v_k\}=\{[\alpha, i]_k|\text{ }\alpha=1,2,\ldots, deg_T([i]_k)-1\}$. Also, we denote the set of \textbf{N}on-\textbf{L}eaf neighbours of $v_k$ to be $NL_k:=\{[i]_k\in V(T)|\ 1\le i\le deg_T(v_k)-1 \text{ and } [i]_k \text{ is not an end-vertex} \}$ for $k=1,2$. Of course, $|NL_k|\ge 1$, where $k=1,2$, for any tree $T$ with $d(T)=5$.

\begin{eg}\label{egA3.1.4}
Let $T$ be the tree of diameter $5$ shown in Figure \ref{figA3.1.1}. We label the vertices as described above. Note also $NL_1=\{[1]_1, [2]_1\}$ and $NL_2=\{[2]_2 \}$.

\begin{center}
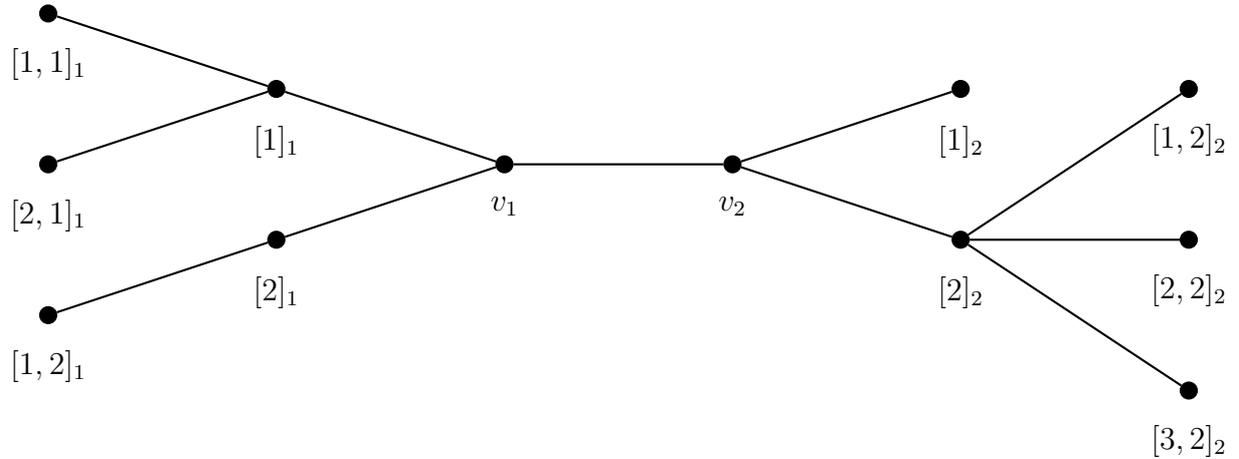

\tikzstyle{every node}=[circle, draw, fill=black!100,
                       inner sep=0pt, minimum width=6pt]
\begin{tikzpicture}[thick,scale=1]%

\draw(-6,2)node[label={[yshift=0cm] 270:{$[1,1]_{1}$}}](11_u1){};
\draw(-6,0)node[label={[yshift=0cm] 270:{$[2,1]_{1}$}}](12_u1){};
\draw(-6,-2)node[label={[yshift=0cm] 270:{$[1,2]_{1}$}}](21_u1){};

\draw(-3,1)node[label={[yshift=-0.2cm] 270:{$[1]_{1}$}}](1_u1){};
\draw(-3,-1)node[label={[yshift=-0.2cm] 270:{$[2]_{1}$}}](2_u1){};
\draw(0,0)node[label={[yshift=-0.2cm]270:{$v_1$}}](u1){};

\draw(3,0)node[label={[yshift=-0.2cm] 270:{$v_2$}}](u2){};
\draw(6,1)node[label={[yshift=-0.2cm] 270:{$[1]_{2}$}}](1_u2){};

\draw(6,-1)node[label={[yshift=-0.2cm] 270:{$[2]_{2}$}}](2_u2){};
\draw(9,1)node[label={[yshift=0cm] 270:{$[1,2]_{2}$}}](12_u2){};
\draw(9,-1)node[label={[yshift=0cm] 270:{$[2,2]_{2}$}}](22_u2){};
\draw(9,-3)node[label={[yshift=0cm] 270:{$[3,2]_{2}$}}](32_u2){};

\draw(u1)--(u2);

\draw(u1)--(1_u1);
\draw(u1)--(2_u1);
\draw(1_u1)--(11_u1);
\draw(1_u1)--(12_u1);
\draw(2_u1)--(21_u1);

\draw(u2)--(1_u2);
\draw(u2)--(2_u2);
\draw(2_u2)--(12_u2);
\draw(2_u2)--(22_u2);
\draw(2_u2)--(32_u2);
\end{tikzpicture}
{\captionof{figure}{Labelling vertices in $T$}\label{figA3.1.1}}
\end{center}
\end{eg}
\indent\par In the following theorem, Koh and Tay \cite{KKM TEG 11} proved some necessary and sufficient conditions for $T(s_1,s_2,\ldots, s_n)\in \mathscr{C}_0$, where $T$ is a tree of diameter $5$.
\begin{thm} (Koh and Tay \cite{KKM TEG 11}) 
\\Let $T$ be a tree with diameter $5$ and its central vertices be $v_1$ and $v_2$. Suppose $A=\{x\in V(T)|\ d_T(x,u)=5=d_T(x,v) \text{ for some }u,v\in V(T), u\neq v\}$. Then, 
\\(a) $T(s_1, s_2,\ldots, s_n)\in \mathscr{C}_0 \cup \mathscr{C}_1$.
\\(b) If $|A|\le 1$, then $T(s_1, s_2,\ldots, s_n) \in \mathscr{C}_0$.
\\(c) If $deg_T(v)\le 2$ for all $v\not\in \{v_1,v_2\}$ and $|A|\ge 2$, then $T^{(2)}\in \mathscr{C}_1$.
\end{thm}
\indent\par In this paper, we give a complete characterisation of vertex-multiplications of trees with diameter $5$ in $\mathscr{C}_0$ and $\mathscr{C}_1$ (see Theorem \ref{thmA3.2.6}).
\section{Main results}
\indent\par For convenience, we shall introduce some notations. Let $D$ be an orientation of $T(s_1,s_2,\ldots,s_n)$ with $s_i\ge2$ for $1\le i\le n$. If $v_p$ and $v_q$, $1\le p, q\le n$ and $p\neq q$, are adjacent vertices in $T$, then for each $i$, $1\le i \le s_p$, we denote by $O_D^{v_q}((i,v_p)):=\{(j,v_q)|\text{ }(i,v_p) \rightarrow (j,v_q), 1\le j \le s_q\}$ and $I_D^{v_q}((i,v_p)):=\{(j,v_q)|\text{ }(j,v_q)\rightarrow (i,v_p), 1\le j \le s_q\}$. If there is no ambiguity, we shall omit the subscript $D$ for the above notations.
\indent\par In the vertex-multiplication graph $G:=T(s_1, s_2,\dots, s_n)$ of $T$, the integer $s_i$ corresponds to the vertex $v_i$, $i=1,2,\ldots,n$. We will loosely use the two denotations of a vertex, for example, if $v_i=[j]_k$, then $s_i=s_{[j]_k}$. For any $v\in V(G)$, we set $(\mathbb{N}_s,v):=\{(1,v),(2,v)\ldots, (s,v)\}$.
\begin{ppn}\label{ppnA3.2.1}
Let $T$ be a tree of diameter $5$. If $s_{1}\ge 3$ or $s_{2}\ge 3$, then $G:=T(s_1,s_2\ldots, s_n)\in \mathscr{C}_0$.
\end{ppn}
\textit{Proof}: WLOG, assume $s_{1}\ge 3$. Let $H$ be a subgraph of $G$, where $s_{1}=3$ and $s_i=2$ for all $i\neq 1$ in $H$. Define an orientation $D$ of $H$ as follows (see Figure \ref{fig:0.2.2}).
\begin{align*}
&(2,[i]_{1})\rightarrow \{(1,[\alpha,i]_{1}), (2,[\alpha,i]_{1})\} \rightarrow (1,[i]_{1}), \\
&(2,v_1)\rightarrow (1,[i]_{1}) \rightarrow \{(1,v_1),(3,v_1)\},\text{ and}\\
&\{(1,v_1),(2,v_1)\}\rightarrow (2,[i]_{1}) \rightarrow (3,v_1),
\end{align*}
for all $i=1,2,\ldots,deg_T(v_1)-1$ and all $\alpha=1,2,\ldots, deg_T([i]_{1})-1$.
\begin{align*}
&(2,[j]_{2})\rightarrow \{(1,[\beta,j]_{2}), (2,[\beta,j]_{2})\} \rightarrow (1,[j]_{2}),\text{ and}\\
&(2,[j]_{2}) \rightarrow (1,v_2)\rightarrow (1,[j]_{2}) \rightarrow (2,v_2) \rightarrow (2,[j]_{2}),
\end{align*}
for all $j=1,2,\ldots,deg_T(v_2)-1$ and all $\beta=1,2,\ldots, deg_T([j]_{2})-1$. Also,
\begin{align*}
&\{(1,v_1),(3,v_1)\}\rightarrow (1,v_2) \rightarrow (2,v_1),\text{ and}\\
&(3,v_1)\rightarrow (2,v_2) \rightarrow \{(1,v_1),(2,v_1)\}.
\end{align*}
\indent\par We emphasize that the above orientation applies to the vertices $[i]_{1}$ and $[j]_{2}$ even if they are end-vertices. It suffices to verify $d(D)=5$ and every vertex in $V(D)$ lies in a directed $C_4$, where $D$ is as shown in Figure \ref{fig:0.2.2}. Then, by Lemma \ref{thmA3.1.3}, $\bar{d}(G)\le \max\{d(D),4\}=5$.
\begin{flushright}
$\Box$
\end{flushright}
\begin{center}
\tikzstyle{every node}=[circle, draw, fill=black!100,
                       inner sep=0pt, minimum width=6pt]
\begin{tikzpicture}[thick,scale=0.75]%
\draw(-6,4)node[label={[yshift=-0.4cm] 90:{$(1,[1,1]_{1})$}}](1_11u1){};
\draw(-6,2)node[label={[yshift=0.4cm] 270:{$(2,[1,1]_{1})$}}](2_11u1){};

\draw(-3,4)node[label={[yshift=0.2cm] 270:{$(1,[1]_{1})$}}](1_1u1){};
\draw(-3,2)node[label={[yshift=0.2cm] 270:{$(2,[1]_{1})$}}](2_1u1){};

\draw(-6,0)node[label={[yshift=0.4cm] 270:{$(1,[1,2]_{1})$}}](1_12u1){};
\draw(-6,-2)node[label={[yshift=0.4cm] 270:{$(2,[1,2]_{1})$}}](2_12u1){};

\draw(-3,0)node[label={[yshift=0.2cm] 270:{$(1,[2]_{1})$}}](1_2u1){};
\draw(-3,-2)node[label={[yshift=0.2cm] 270:{$(2,[2]_{1})$}}](2_2u1){};

\draw(0,3)node[label={[yshift=0cm]90:{$(1,v_1)$}}](1_u1){};
\draw(0,1)node[label={[yshift=0cm]270:{$(2,v_1)$}}](2_u1){};
\draw(0,-1)node[label={[yshift=0cm]270:{$(3,v_1)$}}](3_u1){};

\draw(3,2)node[label={[yshift=0cm] 90:{$(1,v_2)$}}](1_u2){};
\draw(3,0)node[label={[yshift=0cm] 270:{$(2,v_2)$}}](2_u2){};

\draw(6,4)node[label={[yshift=0.2cm] 270:{$(1,[1]_{2})$}}](1_1u2){};
\draw(6,2)node[label={[yshift=0.2cm] 270:{$(2,[1]_{2})$}}](2_1u2){};

\draw(9,4)node[label={[yshift=-0.4cm] 90:{$(1,[1,1]_{2})$}}](1_11u2){};
\draw(9,2)node[label={[yshift=0.4cm] 270:{$(2,[1,1]_{2})$}}](2_11u2){};

\draw(6,0)node[label={[yshift=0.2cm] 270:{$(1,[2]_{2})$}}](1_2u2){};
\draw(6,-2)node[label={[yshift=0.2cm] 270:{$(2,[2]_{2})$}}](2_2u2){};

\draw(9,0)node[label={[yshift=0.4cm] 270:{$(1,[1,2]_{2})$}}](1_12u2){};
\draw(9,-2)node[label={[yshift=0.4cm] 270:{$(2,[1,2]_{2})$}}](2_12u2){};
\draw(-3,8)node[label={[yshift=-0.2cm] 90:{$(1,[3]_{1})$}}](1_3u1){};
\draw(-3,6)node[label={[yshift=0.2cm] 270:{$(2,[3]_{1})$}}](2_3u1){};

\draw(-3,-4)node[label={[yshift=0.2cm] 270:{$(1,[4]_{1})$}}](1_4u1){};
\draw(-3,-6)node[label={[yshift=0.2cm] 270:{$(2,[4]_{1})$}}](2_4u1){};
\draw(6,8)node[label={[yshift=-0.2cm] 90:{$(1,[3]_{2})$}}](1_3u2){};
\draw(6,6)node[label={[yshift=0.2cm] 270:{$(2,[3]_{2})$}}](2_3u2){};

\draw(6,-4)node[label={[yshift=0.2cm] 270:{$(1,[4]_{2})$}}](1_4u2){};
\draw(6,-6)node[label={[yshift=0.2cm] 270:{$(2,[4]_{2})$}}](2_4u2){};

\draw[->, line width=0.3mm, >=latex, shorten <= 0.2cm, shorten >= 0.15cm](1_11u1)--(1_1u1);
\draw[dashed,->, line width=0.3mm, >=latex, shorten <= 0.2cm, shorten >= 0.15cm](2_1u1)--(1_11u1);
\draw[->, line width=0.3mm, >=latex, shorten <= 0.2cm, shorten >= 0.15cm](2_11u1)--(1_1u1);
\draw[dashed,->, line width=0.3mm, >=latex, shorten <= 0.2cm, shorten >= 0.15cm](2_1u1)--(2_11u1);

\draw[->, line width=0.3mm, >=latex, shorten <= 0.2cm, shorten >= 0.15cm](1_1u1)--(1_u1);
\draw[dashed,->, line width=0.3mm, >=latex, shorten <= 0.2cm, shorten >= 0.15cm](2_u1)--(1_1u1);
\draw[->, line width=0.3mm, >=latex, shorten <= 0.2cm, shorten >= 0.15cm](1_1u1)--(3_u1);

\draw[dashed,->, line width=0.3mm, >=latex, shorten <= 0.2cm, shorten >= 0.15cm](1_u1)--(2_1u1);
\draw[dashed,->, line width=0.3mm, >=latex, shorten <= 0.2cm, shorten >= 0.15cm](2_u1)--(2_1u1);
\draw[->, line width=0.3mm, >=latex, shorten <= 0.2cm, shorten >= 0.15cm](2_1u1)--(3_u1);

\draw[->, line width=0.3mm, >=latex, shorten <= 0.2cm, shorten >= 0.15cm](1_3u1)--(1_u1);
\draw[dashed,->, line width=0.3mm, >=latex, shorten <= 0.2cm, shorten >= 0.15cm](2_u1)--(1_3u1);
\draw[->, line width=0.3mm, >=latex, shorten <= 0.2cm, shorten >= 0.15cm](1_3u1)--(3_u1);

\draw[dashed,->, line width=0.3mm, >=latex, shorten <= 0.2cm, shorten >= 0.15cm](1_u1)--(2_3u1);
\draw[dashed,->, line width=0.3mm, >=latex, shorten <= 0.2cm, shorten >= 0.15cm](2_u1)--(2_3u1);
\draw[->, line width=0.3mm, >=latex, shorten <= 0.2cm, shorten >= 0.15cm](2_3u1)--(3_u1);

\draw[->, line width=0.3mm, >=latex, shorten <= 0.2cm, shorten >= 0.15cm](1_4u1)--(1_u1);
\draw[dashed,->, line width=0.3mm, >=latex, shorten <= 0.2cm, shorten >= 0.15cm](2_u1)--(1_4u1);
\draw[->, line width=0.3mm, >=latex, shorten <= 0.2cm, shorten >= 0.15cm](1_4u1)--(3_u1);

\draw[dashed,->, line width=0.3mm, >=latex, shorten <= 0.2cm, shorten >= 0.15cm](1_u1)--(2_4u1);
\draw[dashed,->, line width=0.3mm, >=latex, shorten <= 0.2cm, shorten >= 0.15cm](2_u1)--(2_4u1);
\draw[->, line width=0.3mm, >=latex, shorten <= 0.2cm, shorten >= 0.15cm](2_4u1)--(3_u1);

\draw[dashed,->, line width=0.3mm, >=latex, shorten <= 0.2cm, shorten >= 0.15cm](1_u2)--(2_u1);
\draw[->, line width=0.3mm, >=latex, shorten <= 0.2cm, shorten >= 0.15cm](1_u1)--(1_u2);
\draw[->, line width=0.3mm, >=latex, shorten <= 0.2cm, shorten >= 0.15cm](3_u1)--(1_u2);

\draw[dashed,->, line width=0.3mm, >=latex, shorten <= 0.2cm, shorten >= 0.15cm](2_u2)--(1_u1);
\draw[dashed,->, line width=0.3mm, >=latex, shorten <= 0.2cm, shorten >= 0.15cm](2_u2)--(2_u1);
\draw[->, line width=0.3mm, >=latex, shorten <= 0.2cm, shorten >= 0.15cm](3_u1)--(2_u2);

\draw[->, line width=0.3mm, >=latex, shorten <= 0.2cm, shorten >= 0.15cm](1_u2)--(1_1u2);
\draw[->, line width=0.3mm, >=latex, shorten <= 0.2cm, shorten >= 0.15cm](2_u2)--(2_1u2);
\draw[dashed,->, line width=0.3mm, >=latex, shorten <= 0.2cm, shorten >= 0.15cm](1_1u2)--(2_u2);
\draw[dashed,->, line width=0.3mm, >=latex, shorten <= 0.2cm, shorten >= 0.15cm](2_1u2)--(1_u2);

\draw[->, line width=0.3mm, >=latex, shorten <= 0.2cm, shorten >= 0.15cm](2_1u2)--(1_11u2);
\draw[dashed,->, line width=0.3mm, >=latex, shorten <= 0.2cm, shorten >= 0.15cm](1_11u2)--(1_1u2);
\draw[->, line width=0.3mm, >=latex, shorten <= 0.2cm, shorten >= 0.15cm](2_1u2)--(2_11u2);
\draw[dashed,->, line width=0.3mm, >=latex, shorten <= 0.2cm, shorten >= 0.15cm](2_11u2)--(1_1u2);
\draw[->, line width=0.3mm, >=latex, shorten <= 0.2cm, shorten >= 0.15cm](1_12u1)--(1_2u1);
\draw[->, line width=0.3mm, >=latex, shorten <= 0.2cm, shorten >= 0.15cm](2_12u1)--(1_2u1);
\draw[dashed,->, line width=0.3mm, >=latex, shorten <= 0.2cm, shorten >= 0.15cm](2_2u1)--(1_12u1);
\draw[dashed,->, line width=0.3mm, >=latex, shorten <= 0.2cm, shorten >= 0.15cm](2_2u1)--(2_12u1);

\draw[->, line width=0.3mm, >=latex, shorten <= 0.2cm, shorten >= 0.15cm](1_2u1)--(1_u1);
\draw[dashed,->, line width=0.3mm, >=latex, shorten <= 0.2cm, shorten >= 0.15cm](2_u1)--(1_2u1);
\draw[->, line width=0.3mm, >=latex, shorten <= 0.2cm, shorten >= 0.15cm](1_2u1)--(3_u1);

\draw[dashed,->, line width=0.3mm, >=latex, shorten <= 0.2cm, shorten >= 0.15cm](1_u1)--(2_2u1);
\draw[dashed,->, line width=0.3mm, >=latex, shorten <= 0.2cm, shorten >= 0.15cm](2_u1)--(2_2u1);
\draw[->, line width=0.3mm, >=latex, shorten <= 0.2cm, shorten >= 0.15cm](2_2u1)--(3_u1);
\draw[->, line width=0.3mm, >=latex, shorten <= 0.2cm, shorten >= 0.15cm](1_u2)--(1_2u2);
\draw[->, line width=0.3mm, >=latex, shorten <= 0.2cm, shorten >= 0.15cm](2_u2)--(2_2u2);
\draw[dashed,->, line width=0.3mm, >=latex, shorten <= 0.2cm, shorten >= 0.15cm](1_2u2)--(2_u2);
\draw[dashed,->, line width=0.3mm, >=latex, shorten <= 0.2cm, shorten >= 0.15cm](2_2u2)--(1_u2);

\draw[->, line width=0.3mm, >=latex, shorten <= 0.2cm, shorten >= 0.15cm](2_2u2)--(1_12u2);
\draw[dashed,->, line width=0.3mm, >=latex, shorten <= 0.2cm, shorten >= 0.15cm](1_12u2)--(1_2u2);
\draw[->, line width=0.3mm, >=latex, shorten <= 0.2cm, shorten >= 0.15cm](2_2u2)--(2_12u2);
\draw[dashed,->, line width=0.3mm, >=latex, shorten <= 0.2cm, shorten >= 0.15cm](2_12u2)--(1_2u2);

\draw[->, line width=0.3mm, >=latex, shorten <= 0.2cm, shorten >= 0.15cm](1_u2)--(1_3u2);
\draw[->, line width=0.3mm, >=latex, shorten <= 0.2cm, shorten >= 0.15cm](2_u2)--(2_3u2);
\draw[dashed,->, line width=0.3mm, >=latex, shorten <= 0.2cm, shorten >= 0.15cm](1_3u2)--(2_u2);
\draw[dashed,->, line width=0.3mm, >=latex, shorten <= 0.2cm, shorten >= 0.15cm](2_3u2)--(1_u2);

\draw[->, line width=0.3mm, >=latex, shorten <= 0.2cm, shorten >= 0.15cm](1_u2)--(1_4u2);
\draw[->, line width=0.3mm, >=latex, shorten <= 0.2cm, shorten >= 0.15cm](2_u2)--(2_4u2);
\draw[dashed,->, line width=0.3mm, >=latex, shorten <= 0.2cm, shorten >= 0.15cm](1_4u2)--(2_u2);
\draw[dashed,->, line width=0.3mm, >=latex, shorten <= 0.2cm, shorten >= 0.15cm](2_4u2)--(1_u2);
\end{tikzpicture}
\captionsetup{justification=centering}
{
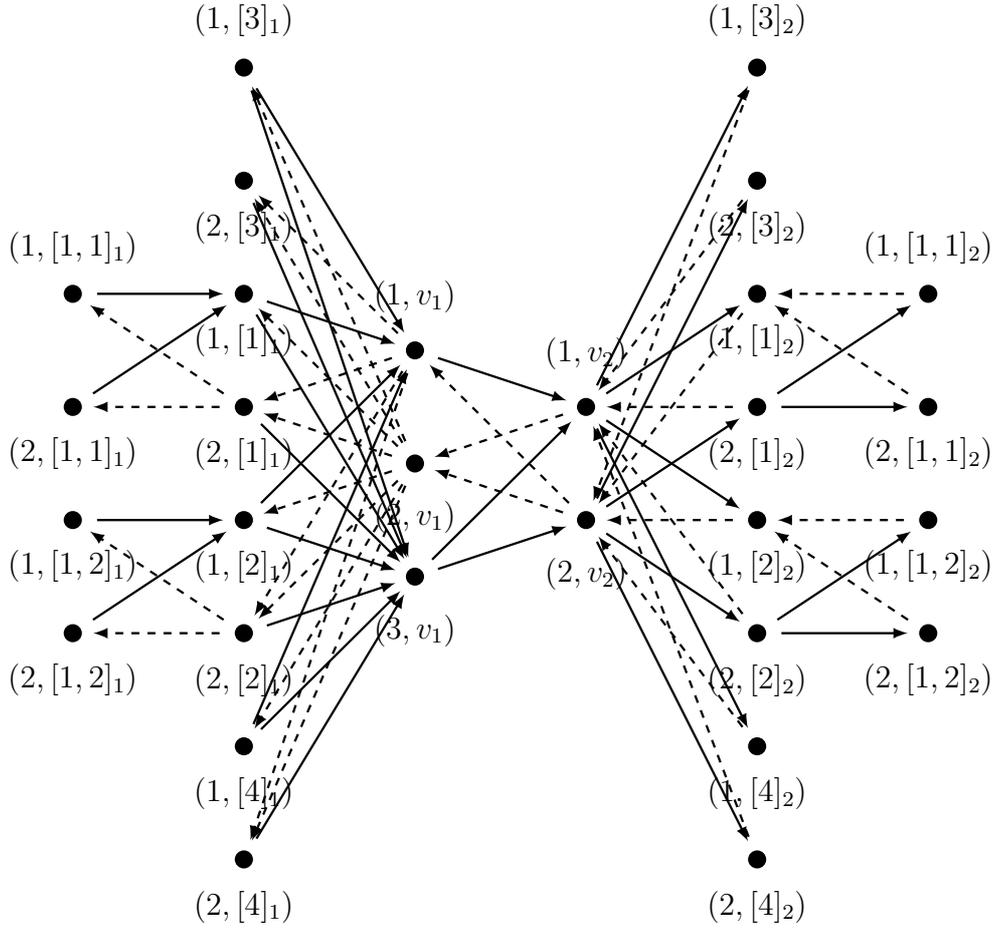
\captionof{figure}{Orientation $D$ for $H$.}\label{fig:0.2.2}}
\end{center}
\indent\par The next two lemmas will be found useful in proving our propositions to come.
\begin{lem}\label{lemA3.2.2}
Let $D$ be an orientation of $T(s_1,s_2,\ldots, s_n)$ with $s_{1}=s_{2}=2$. If $d(D)=5$, then 
\\(i) $|O^{v_1}((p,[i]_{1}))|\ge 1$ and $|I^{v_1}((p,[i]_{1}))|\ge 1$ for all $i=1,2,\ldots, deg_T(v_1)-1$ and $p=1,2,\ldots, s_{[i]_{1}}$, and
\\(ii) $|O^{v_2}((q,[j]_{2}))|\ge 1$ and $|I^{v_2}((q,[j]_{2}))|\ge 1$ for all $j=1,2,\ldots, deg_T(v_2)-1$ and $q=1,2,\ldots, s_{[j]_{2}}$.
\end{lem}
\textit{Proof}:
(i) follows from the fact that $d_D((p,[i]_{1}),(1,[1,j]_{2}))=d_D((1,[1,j]_{2}),(p,[1]_{1}))=4$, where $[j]_2$ is not an end-vertex in $T$. By symmetry, (ii) follows from (i).
\begin{flushright}
$\Box$
\end{flushright}

\begin{lem}\textbf{(Duality)}
\label{lem2.2.1}
~\\Let $D$ be an orientation of a graph $G$. Let $\tilde{D}$ be the orientation of $G$ such that $uv\in E(\tilde{D})$ if and only if $vu \in E(D)$. Then, $d(\tilde{D})=d(D)$.
\end{lem}
\textit{Proof}: Suppose not. Then, there exist some vertices $u,v\in V(\tilde{D})$ such that $d_{\tilde{D}}(u,v)>d(D)$. By definition of $\tilde{D}$, $d_D(v,u)=d_{\tilde{D}}(u,v)$. It follows that $d_{D}(v,u)>d(D)$, a contradiction.
\begin{flushright}
$\Box$
\end{flushright}

\begin{ppn}\label{ppnA3.2.4}
Let $T$ be a tree of diameter $5$ and $m_k=\min\{s_{[i]_{k}}|\ [i]_k\in NL_k\}$ for $k=1,2$. Suppose $G:=T(s_1,s_2,\ldots, s_n)$ satisfy $2\le m_{1}\le 3$ or $2\le m_{2}\le 3$, and $s_1=s_2=2$. Then, 
\begin{align*}
G\in \mathscr{C}_0 \iff  |NL_k|=1 \text{ for some } k=1,2.
\end{align*}
\end{ppn}
\textit{Proof}: ($\Rightarrow$) Suppose $|NL_k|\ge 2 $ for all $k=1,2$. WLOG, assume $2\le m_1\le 3$, and $m_k=s_{[1]_k}$ where $[1]_k\in NL_k$ for $k=1,2$.
\\
\\Case 1. $|I^{v_2}((p,v_1))|=2$ for some $p=1,2$.
\indent\par WLOG, assume $\{(1,v_2),(2,v_2)\} \rightarrow (2,v_1)$. Since for all $i=1,2,\ldots, deg_T(v_1)-1$, and $q=1,2,\ldots, s_{[i]_{1}}$, $d_D((q,[i]_{1}),(1,[1,1]_{2}))\le 5$ and by Lemma \ref{lemA3.2.2}, we have $(2,v_1)\rightarrow (q,[i]_{1})\rightarrow (1,v_1)$. Then, $d_D((1,[1,1]_{1}),(1,[1,i]_{1}))>5$ for any $[i]_1\in NL_1$, $i\neq 1$, a contradiction.
\\
\\Case 2. $|O^{v_2}((p,v_1))|=2$ for some $p=1,2$.
\indent\par By Duality Lemma and Case 1, this case follows.
\\
\\Case 3. $|I^{v_1}((p,v_2))|=|I^{v_2}((p,v_1))|=1$ for all $p=1,2$.
\indent\par WLOG, assume $(1,v_1)\rightarrow (1,v_2)\rightarrow (2,v_1)\rightarrow (2,v_2)\rightarrow (1,v_1)$. 
\\
\\Subcase 3.1. $|O((1,[1,1]_{1}))|=1$.
\indent\par Assume WLOG that $O((1,[1,1]_{1}))=\{(1,[1]_{1})\}$. By Lemma \ref{lemA3.2.2}, $|O^{v_1}((1,[1]_{1}))|=1$. By symmetry, we may assume $O^{v_1}((1,[1]_{1}))=\{(1,v_1)\}$. Since for all $i=1,2,\ldots, deg_T(v_2)-1$, and $q=1,2,\ldots, s_{[i]_{2}}$,  $d_D((1,[1,1]_{1}),(q,[i]_{2}))\le 5$ and by Lemma \ref{lemA3.2.2}, we have $(1,v_2)\rightarrow (q,[i]_{2}) \rightarrow (2,v_2)$. However, $d_D((1,[1,1]_{2}),(1,[1,j]_{2}))>5$ for any $[j]_2\in NL_2$, $j\neq 1$, a contradiction.
\\
\\Subcase 3.2. $|I((1,[1,1]_{1}))|=1$.
\indent\par By Duality Lemma and Subcase 3.1, this subcase also results in a contradiction.

\indent\par ($\Leftarrow$) Assume WLOG that $|NL_1|=1$. Let $H$ be a subgraph of $G$, where $s_i=2$ for all $i$ in $H$. Define an orientation $D$ of $H$ as follows (see Figure \ref{figA3.2.3}).
\begin{align*}
&(1,[\alpha,1]_{1})\rightarrow (1,[1]_{1})\rightarrow (2,[\alpha,1]_{1}) \rightarrow (2,[1]_{1})\rightarrow (1,[\alpha,1]_{1}),\text{ and}\\
&\{(1,[i]_{1}),(2,[i]_{1})\}\rightarrow (1,v_1)\rightarrow \{(1,v_2),(2,v_2)\}\rightarrow (2,v_1)\rightarrow \{(1,[i]_{1}),(2,[i]_{1})\},\\
\end{align*}
for all $i=1,2,\ldots, deg_T(v_1)-1$ and $\alpha=1,2,\ldots, deg_T([1]_{1})-1$.

\begin{align*}
&(2,[j]_{2})\rightarrow \{(1,[\beta, j]_{2}), (2,[\beta, j]_{2})\}\rightarrow (1,[j]_{2}),\text{ and}\\
&(1,[j]_{2}) \rightarrow (1,v_2)\rightarrow (2,[j]_{2}) \rightarrow (2,v_2) \rightarrow (1,[j]_{2}),
\end{align*}
for all $j=1,2,\ldots, deg_T(v_2)-1$ and $\beta=1,2,\ldots, deg_T([j]_{2})-1$.
\indent\par We emphasize that the above orientation applies to the vertices $[i]_{1}$ and $[j]_{2}$ even if they are end-vertices. It suffices to verify $d(D)=5$ and every vertex in $V(D)$ lies in a directed $C_4$, where $D$ is as shown in Figure \ref{figA3.2.3}. Then, by Lemma \ref{thmA3.1.3}, $\bar{d}(G)\le \max\{d(D),4\}=5$.
\begin{flushright}
$\Box$
\end{flushright}

\begin{center}
\tikzstyle{every node}=[circle, draw, fill=black!100,
                       inner sep=0pt, minimum width=6pt]
\begin{tikzpicture}[thick,scale=0.7]%
\draw(-6,2)node[label={[yshift=-0.4cm] 90:{$(1,[1,1]_{1})$}}](1_11u1){};
\draw(-6,0)node[label={[yshift=0.4cm] 270:{$(2,[1,1]_{1})$}}](2_11u1){};

\draw(-3,2)node[label={[yshift=0.2cm] 270:{$(1,[1]_{1})$}}](1_1u1){};
\draw(-3,0)node[label={[yshift=0.2cm] 270:{$(2,[1]_{1})$}}](2_1u1){};

\draw(0,2)node[label={[yshift=0cm]90:{$(1,v_1)$}}](1_u1){};
\draw(0,0)node[label={[yshift=0.15cm]270:{$(2,v_1)$}}](2_u1){};

\draw(3,2)node[label={[yshift=0cm] 90:{$(1,v_2)$}}](1_u2){};
\draw(3,0)node[label={[yshift=0.15cm] 270:{$(2,v_2)$}}](2_u2){};

\draw(6,4)node[label={[yshift=0.2cm] 270:{$(1,[1]_{2})$}}](1_1u2){};
\draw(6,2)node[label={[yshift=0.2cm] 270:{$(2,[1]_{2})$}}](2_1u2){};

\draw(9,4)node[label={[yshift=-0.4cm] 90:{$(1,[1,1]_{2})$}}](1_11u2){};
\draw(9,2)node[label={[yshift=0.4cm] 270:{$(2,[1,1]_{2})$}}](2_11u2){};

\draw(6,0)node[label={[yshift=0.2cm] 270:{$(1,[2]_{2})$}}](1_2u2){};
\draw(6,-2)node[label={[yshift=0.2cm] 270:{$(2,[2]_{2})$}}](2_2u2){};

\draw(9,0)node[label={[yshift=0.4cm] 270:{$(1,[1,2]_{2})$}}](1_12u2){};
\draw(9,-2)node[label={[yshift=0.4cm] 270:{$(2,[1,2]_{2})$}}](2_12u2){};

\draw(-3,6)node[label={[yshift=-0.2cm] 90:{$(1,[2]_{1})$}}](1_2u1){};
\draw(-3,4)node[label={[yshift=0.2cm] 270:{$(2,[2]_{1})$}}](2_2u1){};

\draw(-3,-2)node[label={[yshift=0.2cm] 270:{$(1,[3]_{1})$}}](1_3u1){};
\draw(-3,-4)node[label={[yshift=0.2cm] 270:{$(2,[3]_{1})$}}](2_3u1){};

\draw(6,8)node[label={[yshift=-0.2cm] 90:{$(1,[3]_{2})$}}](1_3u2){};
\draw(6,6)node[label={[yshift=0.2cm] 270:{$(2,[3]_{2})$}}](2_3u2){};

\draw(6,-4)node[label={[yshift=0.2cm] 270:{$(1,[4]_{2})$}}](1_4u2){};
\draw(6,-6)node[label={[yshift=0.2cm] 270:{$(2,[4]_{2})$}}](2_4u2){};

\draw[->, line width=0.3mm, >=latex, shorten <= 0.2cm, shorten >= 0.15cm](1_11u1)--(1_1u1);
\draw[dashed,->, line width=0.3mm, >=latex, shorten <= 0.2cm, shorten >= 0.15cm](2_1u1)--(1_11u1);
\draw[->, line width=0.3mm, >=latex, shorten <= 0.2cm, shorten >= 0.15cm](2_11u1)--(2_1u1);
\draw[dashed,->, line width=0.3mm, >=latex, shorten <= 0.2cm, shorten >= 0.15cm](1_1u1)--(2_11u1);

\draw[->, line width=0.3mm, >=latex, shorten <= 0.2cm, shorten >= 0.15cm](1_1u1)--(1_u1);
\draw[dashed,->, line width=0.3mm, >=latex, shorten <= 0.2cm, shorten >= 0.15cm](2_u1)--(1_1u1);
\draw[->, line width=0.3mm, >=latex, shorten <= 0.2cm, shorten >= 0.15cm](2_1u1)--(1_u1);
\draw[dashed,->, line width=0.3mm, >=latex, shorten <= 0.2cm, shorten >= 0.15cm](2_u1)--(2_1u1);

\draw[->, line width=0.3mm, >=latex, shorten <= 0.2cm, shorten >= 0.15cm](1_u1)--(1_u2);
\draw[->, line width=0.3mm, >=latex, shorten <= 0.2cm, shorten >= 0.15cm](1_u1)--(2_u2);
\draw[dashed, ->, line width=0.3mm, >=latex, shorten <= 0.2cm, shorten >= 0.15cm](1_u2)--(2_u1);
\draw[dashed, ->, line width=0.3mm, >=latex, shorten <= 0.2cm, shorten >= 0.15cm](2_u2)--(2_u1);

\draw[dashed,->, line width=0.3mm, >=latex, shorten <= 0.2cm, shorten >= 0.15cm](1_1u2)--(1_u2);
\draw[->, line width=0.3mm, >=latex, shorten <= 0.2cm, shorten >= 0.15cm](1_u2)--(2_1u2);
\draw[dashed,->, line width=0.3mm, >=latex, shorten <= 0.2cm, shorten >= 0.15cm](2_1u2)--(2_u2);
\draw[->, line width=0.3mm, >=latex, shorten <= 0.2cm, shorten >= 0.15cm](2_u2)--(1_1u2);

\draw[dashed,->, line width=0.3mm, >=latex, shorten <= 0.2cm, shorten >= 0.15cm](1_11u2)--(1_1u2);
\draw[dashed,->, line width=0.3mm, >=latex, shorten <= 0.2cm, shorten >= 0.15cm](2_11u2)--(1_1u2);
\draw[->, line width=0.3mm, >=latex, shorten <= 0.2cm, shorten >= 0.15cm](2_1u2)--(1_11u2);
\draw[->, line width=0.3mm, >=latex, shorten <= 0.2cm, shorten >= 0.15cm](2_1u2)--(2_11u2);
\draw[dashed,->, line width=0.3mm, >=latex, shorten <= 0.2cm, shorten >= 0.15cm](1_2u2)--(1_u2);
\draw[->, line width=0.3mm, >=latex, shorten <= 0.2cm, shorten >= 0.15cm](1_u2)--(2_2u2);
\draw[dashed,->, line width=0.3mm, >=latex, shorten <= 0.2cm, shorten >= 0.15cm](2_2u2)--(2_u2);
\draw[->, line width=0.3mm, >=latex, shorten <= 0.2cm, shorten >= 0.15cm](2_u2)--(1_2u2);

\draw[dashed,->, line width=0.3mm, >=latex, shorten <= 0.2cm, shorten >= 0.15cm](1_12u2)--(1_2u2);
\draw[dashed,->, line width=0.3mm, >=latex, shorten <= 0.2cm, shorten >= 0.15cm](2_12u2)--(1_2u2);
\draw[->, line width=0.3mm, >=latex, shorten <= 0.2cm, shorten >= 0.15cm](2_2u2)--(1_12u2);
\draw[->, line width=0.3mm, >=latex, shorten <= 0.2cm, shorten >= 0.15cm](2_2u2)--(2_12u2);

\draw[dashed,->, line width=0.3mm, >=latex, shorten <= 0.2cm, shorten >= 0.15cm](1_3u2)--(1_u2);
\draw[->, line width=0.3mm, >=latex, shorten <= 0.2cm, shorten >= 0.15cm](1_u2)--(2_3u2);
\draw[dashed,->, line width=0.3mm, >=latex, shorten <= 0.2cm, shorten >= 0.15cm](2_3u2)--(2_u2);
\draw[->, line width=0.3mm, >=latex, shorten <= 0.2cm, shorten >= 0.15cm](2_u2)--(1_3u2);

\draw[dashed,->, line width=0.3mm, >=latex, shorten <= 0.2cm, shorten >= 0.15cm](1_4u2)--(1_u2);
\draw[->, line width=0.3mm, >=latex, shorten <= 0.2cm, shorten >= 0.15cm](1_u2)--(2_4u2);
\draw[dashed,->, line width=0.3mm, >=latex, shorten <= 0.2cm, shorten >= 0.15cm](2_4u2)--(2_u2);
\draw[->, line width=0.3mm, >=latex, shorten <= 0.2cm, shorten >= 0.15cm](2_u2)--(1_4u2);

\draw[->, line width=0.3mm, >=latex, shorten <= 0.2cm, shorten >= 0.15cm](1_2u1)--(1_u1);
\draw[dashed,->, line width=0.3mm, >=latex, shorten <= 0.2cm, shorten >= 0.15cm](2_u1)--(1_2u1);
\draw[->, line width=0.3mm, >=latex, shorten <= 0.2cm, shorten >= 0.15cm](2_2u1)--(1_u1);
\draw[dashed,->, line width=0.3mm, >=latex, shorten <= 0.2cm, shorten >= 0.15cm](2_u1)--(2_2u1);

\draw[->, line width=0.3mm, >=latex, shorten <= 0.2cm, shorten >= 0.15cm](1_3u1)--(1_u1);
\draw[dashed,->, line width=0.3mm, >=latex, shorten <= 0.2cm, shorten >= 0.15cm](2_u1)--(1_3u1);
\draw[->, line width=0.3mm, >=latex, shorten <= 0.2cm, shorten >= 0.15cm](2_3u1)--(1_u1);
\draw[dashed,->, line width=0.3mm, >=latex, shorten <= 0.2cm, shorten >= 0.15cm](2_u1)--(2_3u1);

\end{tikzpicture}
\captionsetup{justification=centering}
{
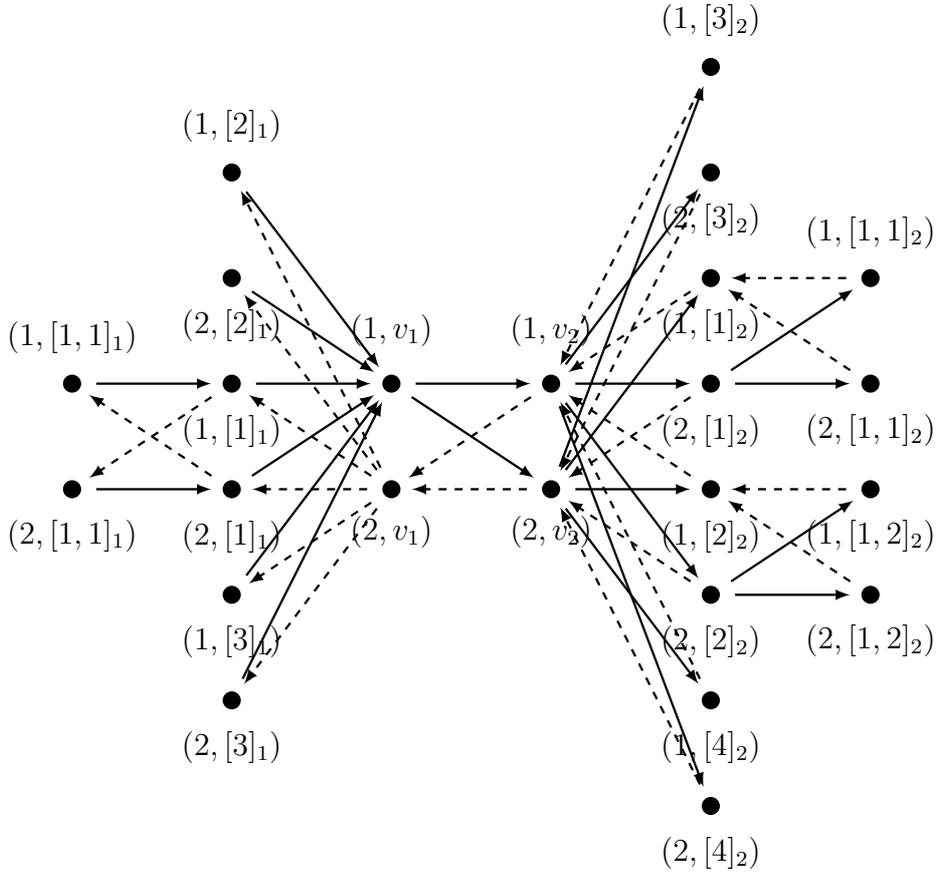
\captionof{figure}{Orientation $D$ for $H$.}\label{figA3.2.3}}
\end{center}

\begin{ppn}\label{ppn6.2.5}
Let $T$ be a tree of diameter $5$ and $m_k=\min\{s_{[i]_{k}}|\ [i]_k\in NL_k\}$ for $k=1,2$. If $m_1\ge 4$ and $m_2\ge 4$, then $G:=T(s_1,s_2\ldots, s_n)\in \mathscr{C}_0$.
\end{ppn}
\textit{Proof}: Let $H$ be a subgraph of $G$, where $s_{[i]_k}=4$ for all $[i]_k\in NL_k$, $k=1,2$, and $s_j=2$ otherwise. Define an orientation $D$ of $H$ as follows (see Figure \ref{figA3.2.4}). For $k=1,2$,
\begin{align*}
&\{(3,[i]_{k}), (4,[i]_{k})\}\rightarrow \{(1,[\alpha,i]_{k}), (2,[\alpha,i]_{k})\} \rightarrow \{(1,[i]_{k}), (2,[i]_{k})\},\\
&(2,v_k)\rightarrow \{(1,[i]_{k}), (3,[i]_{k})\} \rightarrow (1,v_k)\rightarrow \{(2,[i]_{k}), (4,[i]_{k})\}\rightarrow (2,v_k),\text{ and}\\
&(1,v_1) \rightarrow (1,v_2)\rightarrow (2,v_1) \rightarrow (2,v_2) \rightarrow (1,v_1),
\end{align*}
for all $i=1,2,\ldots, deg_T(v_k)-1$ and $\alpha=1,2,\ldots, deg_T([i]_{k})-1$. 

\indent\par We emphasize that the above orientation applies to the vertices $[i]_{1}$ and $[j]_{2}$ even if they are end-vertices. Note that for any $p,q,r,k=1,2$, $1\le i, j\le deg_T(v_k)$, $i\neq j$, $\alpha=1,2,\ldots, deg_T([i]_{k})-1$, and $\beta=1,2,\ldots, deg_T([j]_{k})-1$, we have $d_D((p,[\alpha,i]_k), (q,v_k))=2$ and $d_D((q,v_k),(r,[\beta,j]_k))=2$ so that $d_D((p,[\alpha,i]_k), (r,[\beta,j]_k))\le d_D((p,[\alpha,i]_k), (q,v_k))+d_D((q,v_k),(r,[\beta,j]_k))=4$.
Furthermore, in view of the orientation's similarity, it suffices to verify $d(D)=5$ and every vertex in $V(D)$ lies in a directed $C_4$, where $D$ is as shown in Figure \ref{figA3.2.4}. Hence, by Lemma \ref{thmA3.1.3}, $\bar{d}(G)\le \max\{d(D),4\}=5$.
\begin{flushright}
$\Box$
\end{flushright}

\begin{center}
\tikzstyle{every node}=[circle, draw, fill=black!100,
                       inner sep=0pt, minimum width=6pt]
\begin{tikzpicture}[thick,scale=1]%
\draw(-6,2)node[label={[yshift=-0.4cm] 90:{$(1,[1,1]_{1})$}}](1_11u1){};
\draw(-6,0)node[label={[yshift=0.4cm] 270:{$(2,[1,1]_{1})$}}](2_11u1){};

\draw(-3,4)node[label={[yshift=-0.2cm] 90:{$(1,[1]_{1})$}}](1_1u1){};
\draw(-3,2)node[label={[yshift=-0.2cm] 90:{$(2,[1]_{1})$}}](2_1u1){};
\draw(-3,0)node[label={[yshift=0.2cm] 270:{$(3,[1]_{1})$}}](3_1u1){};
\draw(-3,-2)node[label={[yshift=0.2cm] 270:{$(4,[1]_{1})$}}](4_1u1){};

\draw(0,2)node[label={[yshift=-0.1cm]90:{$(1,v_1)$}}](1_u1){};
\draw(0,0)node[label={[yshift=0.15cm]270:{$(2,v_1)$}}](2_u1){};

\draw(3,2)node[label={[yshift=-0.1cm] 90:{$(1,v_2)$}}](1_u2){};
\draw(3,0)node[label={[yshift=0.15cm] 270:{$(2,v_2)$}}](2_u2){};

\draw(6,4)node[label={[yshift=-0.2cm] 90:{$(1,[1]_{2})$}}](1_1u2){};
\draw(6,2)node[label={[yshift=-0.2cm] 90:{$(2,[1]_{2})$}}](2_1u2){};
\draw(6,0)node[label={[yshift=0.2cm] 270:{$(3,[1]_{2})$}}](3_1u2){};
\draw(6,-2)node[label={[yshift=0.2cm] 270:{$(4,[1]_{2})$}}](4_1u2){};

\draw(9,2)node[label={[yshift=-0.4cm] 90:{$(1,[1,1]_{2})$}}](1_11u2){};
\draw(9,0)node[label={[yshift=0.4cm] 270:{$(2,[1,1]_{2})$}}](2_11u2){};

\draw[->, line width=0.3mm, >=latex, shorten <= 0.2cm, shorten >= 0.15cm](1_11u1)--(1_1u1);
\draw[->, line width=0.3mm, >=latex, shorten <= 0.2cm, shorten >= 0.15cm](1_11u1)--(2_1u1);
\draw[dashed,->, line width=0.3mm, >=latex, shorten <= 0.2cm, shorten >= 0.15cm](3_1u1)--(1_11u1);
\draw[dashed,->, line width=0.3mm, >=latex, shorten <= 0.2cm, shorten >= 0.15cm](4_1u1)--(1_11u1);

\draw[->, line width=0.3mm, >=latex, shorten <= 0.2cm, shorten >= 0.15cm](2_11u1)--(1_1u1);
\draw[->, line width=0.3mm, >=latex, shorten <= 0.2cm, shorten >= 0.15cm](2_11u1)--(2_1u1);
\draw[dashed,->, line width=0.3mm, >=latex, shorten <= 0.2cm, shorten >= 0.15cm](3_1u1)--(2_11u1);
\draw[dashed,->, line width=0.3mm, >=latex, shorten <= 0.2cm, shorten >= 0.15cm](4_1u1)--(2_11u1);

\draw[->, line width=0.3mm, >=latex, shorten <= 0.2cm, shorten >= 0.15cm](1_1u1)--(1_u1);
\draw[dashed,->, line width=0.3mm, >=latex, shorten <= 0.2cm, shorten >= 0.15cm](2_u1)--(1_1u1);
\draw[dashed,->, line width=0.3mm, >=latex, shorten <= 0.2cm, shorten >= 0.15cm](1_u1)--(2_1u1);
\draw[->, line width=0.3mm, >=latex, shorten <= 0.2cm, shorten >= 0.15cm](2_1u1)--(2_u1);

\draw[->, line width=0.3mm, >=latex, shorten <= 0.2cm, shorten >= 0.15cm](3_1u1)--(1_u1);
\draw[dashed,->, line width=0.3mm, >=latex, shorten <= 0.2cm, shorten >= 0.15cm](2_u1)--(3_1u1);
\draw[dashed,->, line width=0.3mm, >=latex, shorten <= 0.2cm, shorten >= 0.15cm](1_u1)--(4_1u1);
\draw[->, line width=0.3mm, >=latex, shorten <= 0.2cm, shorten >= 0.15cm](4_1u1)--(2_u1);

\draw[->, line width=0.3mm, >=latex, shorten <= 0.2cm, shorten >= 0.15cm](1_u1)--(1_u2);
\draw[dashed,->, line width=0.3mm, >=latex, shorten <= 0.2cm, shorten >= 0.15cm](2_u2)--(1_u1);
\draw[dashed,->, line width=0.3mm, >=latex, shorten <= 0.2cm, shorten >= 0.15cm](1_u2)--(2_u1);
\draw[->, line width=0.3mm, >=latex, shorten <= 0.2cm, shorten >= 0.15cm](2_u1)--(2_u2);

\draw[dashed,->, line width=0.3mm, >=latex, shorten <= 0.2cm, shorten >= 0.15cm](1_1u2)--(1_u2);
\draw[dashed,->, line width=0.3mm, >=latex, shorten <= 0.2cm, shorten >= 0.15cm](2_1u2)--(2_u2);
\draw[->, line width=0.3mm, >=latex, shorten <= 0.2cm, shorten >= 0.15cm](1_u2)--(2_1u2);
\draw[->, line width=0.3mm, >=latex, shorten <= 0.2cm, shorten >= 0.15cm](2_u2)--(1_1u2);
\draw[dashed,->, line width=0.3mm, >=latex, shorten <= 0.2cm, shorten >= 0.15cm](3_1u2)--(1_u2);
\draw[dashed,->, line width=0.3mm, >=latex, shorten <= 0.2cm, shorten >= 0.15cm](4_1u2)--(2_u2);
\draw[->, line width=0.3mm, >=latex, shorten <= 0.2cm, shorten >= 0.15cm](1_u2)--(4_1u2);
\draw[->, line width=0.3mm, >=latex, shorten <= 0.2cm, shorten >= 0.15cm](2_u2)--(3_1u2);

\draw[dashed,->, line width=0.3mm, >=latex, shorten <= 0.2cm, shorten >= 0.15cm](1_11u2)--(1_1u2);
\draw[dashed,->, line width=0.3mm, >=latex, shorten <= 0.2cm, shorten >= 0.15cm](1_11u2)--(2_1u2);
\draw[->, line width=0.3mm, >=latex, shorten <= 0.2cm, shorten >= 0.15cm](3_1u2)--(1_11u2);
\draw[->, line width=0.3mm, >=latex, shorten <= 0.2cm, shorten >= 0.15cm](4_1u2)--(1_11u2);

\draw[dashed,->, line width=0.3mm, >=latex, shorten <= 0.2cm, shorten >= 0.15cm](2_11u2)--(1_1u2);
\draw[dashed,->, line width=0.3mm, >=latex, shorten <= 0.2cm, shorten >= 0.15cm](2_11u2)--(2_1u2);
\draw[->, line width=0.3mm, >=latex, shorten <= 0.2cm, shorten >= 0.15cm](3_1u2)--(2_11u2);
\draw[->, line width=0.3mm, >=latex, shorten <= 0.2cm, shorten >= 0.15cm](4_1u2)--(2_11u2);
\end{tikzpicture}
\captionsetup{justification=centering}
{\captionof{figure}{Orientation $D$ for $H$.}\label{figA3.2.4}}
\end{center}

\indent\par The findings in this section provide a complete characterisation of vertex-multiplications of trees of diameter $5$ in $\mathscr{C}_0$ and $\mathscr{C}_1$. We summarise them in Theorem \ref{thmA3.2.6}.
\newpage
\begin{thm}\label{thmA3.2.6}
Let $G:=T(s_1,s_2\ldots, s_n)$ with $d(T)=5$, and $m_k=\min\{s_{[i]_{k}}|\ i\in NL_k\}$ for $k=1,2$.
\end{thm}
\begin{spacing}{1.2}
\begin{centering}
\begin{tabular}{| P{4.5cm} | P{5cm}| P{3.5cm}|}
\hline
\textbf{Underlying conditions} & $\bm{G\in\mathscr{C}_0}$ \textbf{or} $\bm{\mathscr{C}_1 ?}$ & \textbf{Reference}\\
\hline
$s_{i}\ge 3$ for some $i=1,2$ & $G\in \mathscr{C}_0$. & Proposition \ref{ppnA3.2.1}\\ 
\hline\hline
$m_1\ge 4$ and $m_2\ge 4$ & $G\in \mathscr{C}_0$. & Proposition \ref{ppn6.2.5}\\ 
\hline\hline
$s_{1}=s_{2}=2$, and $2\le m_i\le 3$ for some $i=1,2$ & $G\in \mathscr{C}_0 \iff$ $|NL_k|=1$ for some $k=1,2$. & Proposition \ref{ppnA3.2.4}\\ 
\hline
\end{tabular}
{\captionof{table}{Summary for $T(s_1, s_2, \ldots, s_n)$ with $d(T)=5$.}}
\end{centering}
\end{spacing}

\indent\par Koh and Tay \cite{KKM TEG 11} also showed that vertex-multiplications of trees with diameter at least 6 belong to $\mathscr{C}_0$. Hence, it remains open to characterise vertex-multiplications of trees, with diameter 3 and 4, belonging to $\mathscr{C}_0$. We conclude the paper by proposing the following problem.

\begin{prob}
For trees $T$ with $d(T)=3$ (and $d(T)=4$ resp.), characterise the tree vertex-multiplications $T(s_1,s_2,\ldots,s_n)$ that belong to $\mathscr{C}_0$.  
\end{prob}
\textbf{Acknowledgement}
\indent\par The first author would like to thank the National Institute of Education, Nanyang Technological University of Singapore, for the generous support of the Nanyang Technological University Research Scholarship.

\end{document}